\documentclass[oneside,centertags,12pt]{amsart}
\usepackage{multirow}
\usepackage{graphicx}
\usepackage{float}
\usepackage{array}
\usepackage{lscape}
\usepackage{amsfonts}
\usepackage{amsmath}
\usepackage{amssymb}
\usepackage{amsthm}
\usepackage{hyperref}
\numberwithin{equation}{section}

\theoremstyle{plain}
\newtheorem{thm}{Theorem}[section]
\newtheorem{prop}[thm]{Proposition}
\newtheorem{proposition}[thm]{Proposition}
\newtheorem{lemma}[thm]{Lemma}

\newtheorem{cor}[thm]{Corollary}

\theoremstyle{definition}

\newtheorem{definition1}[thm]{Definition}

\newtheorem{example}[thm]{Example}

\newtheorem{q}[thm]{Question}

\makeindex

\begin{document}

\title{On the divisibility of binomial coefficients}
\author{S\'{\i}lvia Casacuberta}
\address{Harvard University, Cambridge, MA 02138, USA}
\email{scasacubertapuig@college.harvard.edu}
\date{}

\thanks{}

\begin{abstract}
In Pacific J. Math.\ 292 (2018), 223--238, Shareshian and Woodroofe asked if for every positive integer $n$ there exist primes $p$ and $q$ such that, for all integers $k$ with $1 \leq k \leq n-1$, the binomial coefficient $\binom{n}{k}$ is divisible by at least one of $p$ or~$q$. We give conditions under which a number $n$ has this property and discuss a variant of this problem involving more than two primes. We prove that every positive integer $n$ has infinitely many multiples with this property.
\end{abstract}

\maketitle

\section{Introduction}
Binomial coefficients display interesting divisibility properties. Conditions under which a prime power $p^a$ divides a binomial coefficient $\binom{n}{k}$ are given by Kummer's Theorem \cite{kummer} and also by a generalized form of Lucas' Theorem \cite{glucas,lucas}. 

Still, there are problems involving divisibility of binomial coefficients that remain unsolved. In this article we investigate the following question, which was asked by Shareshian and Woodroofe in~\cite{sw}.

\begin{q}
\label{q}
Is it true that for every positive integer $n$ there exist primes $p$ and $q$ such that, for all integers $k$ with $1 \leq k \leq n-1$, the binomial coefficient $\binom{n}{k}$ is divisible by $p$ or~$q$? 
\end{q}

As in \cite{sw}, we say that $n$ \emph{satisfies Condition $1$} if such primes $p$ and $q$ exist for~$n$.
In this article we discuss sufficient conditions under which an integer $n$ satisfies Condition~1.
In Sections \ref{sec2} and \ref{secq2} we prove a variation of the Sieve Lemma from \cite{sw} and use it to show that $n$ satisfies Condition~1 if certain inequalities hold. 
In Section~\ref{sec5} we infer that every positive integer has infinitely many multiples for which Condition~1 is satisfied.

The collection of numbers for which Condition~1 is not known to hold has asymptotic density~$0$ assuming the truth of Cram\'er's conjecture (as first shown in~\cite{sw}) and includes most \emph{primorials} $p_1p_2\cdots p_i$, where $p_1,\dots,p_i$ are the first $i$ primes, namely those primorials such that $(p_1p_2\cdots p_i)-1$ is not a prime. 

In addition, we introduce the following variant of Condition~1:

\begin{definition1}
\label{def1}
A positive integer $n$ satisfies the \emph{$N$-variation} of Condition~1 if there exist $N$ different primes $p_1,\dots,p_N$ such that if $1 \leq k \leq n-1$ then $\binom{n}{k}$ is divisible by at least one of $p_1,\dots,p_N$.
\end{definition1}

For example, it follows from Kummer's Theorem or from Lucas' Theorem that a positive integer $n$ satisfies the $1$-variation of Condition~$1$ if and only if $n$ is a prime power, and every integer $n$ satisfies the $m$-variation of Condition~1 if $n=p_1^{a_1}\cdots p_m^{a_m}$ where $p_1,\dots,p_m$ are distinct primes.
In Section~\ref{sec4} we discuss upper bounds on $N$ so that a given $n$ satisfies the $N$-variation of Condition~1.

\newpage

\section{An extended sieve lemma}
\label{sec2}

Our results in this section will be based on Lucas' Theorem:

\begin{thm}[Lucas \cite{lucas}]
Let $p$ be a prime and let 
\begin{align*}
n & = n_{r}p^{r} + n_{r-1}p^{r-1} + \cdots + n_{1}p + n_{0} \\
k & = k_{r}p^{r} + k_{r-1}p^{r-1} + \cdots + k_{1}p + k_0
\end{align*}
be base $p$ expansions of two positive integers, where $0\leq n_i<p$ and $0\leq k_i<p$ for all~$i$, and $n_r\ne 0$. Then 
\[
\binom{n}{k}  \equiv \prod_{i=0}^{r} \binom{n_{i}}{k_{i}} \; \text{\rm mod~$p$}.
\]
\end{thm}

By convention, a binomial coefficient $\binom{n_i}{k_i}$ is zero if $n_i < k_i$. Hence, if any of the digits of the base~$p$ expansion of $n$ is~$0$ whereas the corresponding digit in the base~$p$ expansion of $k$ is nonzero, then $\binom{n}{k}$ is divisible by~$p$. 
As a particular case, if a prime power $p^a$ with $a>0$ divides $n$ and does not divide~$k$, then $\binom{n}{k}$ is divisible by~$p$.

Observe that, if $n$ satisfies Condition~1 with two primes $p$ and~$q$, then at least one of these primes has to be a divisor of~$n$, because otherwise $\binom{n}{1}$ would not be divisible by any of them. The next two results are elementary consequences of Lucas' Theorem. 


\begin{prop}
If $n=p^{a} + 1$ with $p$ a prime and $a>0$, then $n$ satisfies Condition~$1$ with $p$ and any prime dividing~$n$.
\end{prop}

\begin{proof}
If $n-1$ is a prime power then the two summands in the left-hand term of the equality
\[
\binom{n-1}{k-1}+\binom{n-1}{k}=\binom{n}{k}
\] 
are divisible by $p$ by Lucas' Theorem if $2 \leq k \leq n-2$, and hence 
$\binom{n}{k}$ is also divisible by~$p$.
When $k=1$ or $k=n-1$, we have that $\binom{n}{k} = n$, so any prime factor of $n$ divides $\binom{n}{k}$.
\end{proof}

\begin{prop}
\label{twoprimes}
If a positive integer $n$ is equal to the product of two prime powers $p_{1}^{a}$ and $p_{2}^{b}$ with $a>0$, $b>0$, and $p_1\ne p_2$, then $n$ satisfies Condition~$1$ with $p_{1}$ and~$p_{2}$. 
\end{prop}

\begin{proof}
The base $p_{1}$ expansion of $n$ ends with $a$ zeroes and the base $p_{2}$ expansion of $n$ ends with $b$ zeroes. Because a positive integer $k$ smaller than $n$ cannot be divisible by both $p_{1}^{a}$ and $p_{2}^{b}$, it is not possible that $k$ ends with $a$ zeroes in base~$p_{1}$ and $b$ zeroes in base~$p_{2}$. Consequently, we can apply Lucas' Theorem modulo $p_{1}$ if $p_{1}^{a}$ does not divide~$k$ or modulo $p_{2}$ if $p_{2}^{b}$ does not divide~$k$.
\end{proof}

Proposition~\ref{twoprimes} generalizes as follows.

\begin{prop}
If $p_1,\dots,p_m$ are distinct primes and $n=p_1^{a_1}\cdots p_m^{a_m}$ with $a_i>0$ for all~$i$, then $n$ satisfies the $m$-variation of Condition~$1$ with $p_1\dots,p_m$.
\end{prop}

\begin{proof}
If $1\leq k\leq n-1$, then the base $p_i$ expansion of $k$ ends with less zeroes than the base $p_i$ expansion of $n$ for at least one prime factor $p_i$ of~$n$. 
\end{proof}

The following result extends \cite[Lemma~4.3]{sw}. It is the starting point of our discussion of Question~\ref{q} in the next sections.

\begin{thm}
\label{MIAU}
Let $n$ be a positive integer and suppose that $p^a$ divides $n$ where $p$ is a prime and $a>0$. Suppose that there is a prime $q$ with $n/(d+1)<q<n/d$, where $d\ge 1$. Then $\binom{n}{k}$ is divisible by $p$ or $q$ except possibly when $k$ is a multiple of $p^a$ belonging to one of the intervals $[cq,\,cq+\beta]$ with $\beta=n-dq$ and $0\le c<(d+1)/2$.
\end{thm}

\begin{proof}
By symmetry, we only need to consider those values of $k$ with $k\le n/2$. Moreover, we may restrict our study further to those values of $k$ that are multiples of~$p^a$, since otherwise $\binom{n}{k}$ is divisible by~$p$.

Since $q<n/d$, the number $\beta=n-dq$ is positive. If $k\le\beta$ then $k$ is in the interval $[0,\beta]$, which is the case $c=0$ in the statement of the theorem.

The assumption that $n/(d+1)<q$ is equivalent to assuming the inequality $n-dq<q$, which implies that the last digit in the base $q$ expansion of $n$ is equal to~$\beta$. Hence, if $\beta<k<q$ then we may infer from Lucas' Theorem that $\binom{n}{k}$ is divisible by~$q$. 

The remaining range of values of $k$ to be considered is $q\le k\le n/2$. In this case we look at the last digit of the base~$q$ expansion of~$k$. If this last digit is bigger than~$\beta$, then $\binom{n}{k}$ is again divisible by~$q$. Thus the undecided cases are those in which the residue of $k$ mod~$q$ is smaller than or equal to~$\beta$. This happens when $cq\le k\le cq+\beta$ for some positive integer~$c$, and if $cq\le k\le n/2$ then $c\le n/(2q)<(d+1)/2$.
\end{proof}

By the Bertrand--Chebyshev Theorem \cite{bertrand}, for every integer $n > 2$ there exists a prime $q$ such that $n/2 < q < n$. This yields the following particular instance of Theorem~\ref{MIAU}, which is also a special case of \cite[Lemma~4.3]{sw}.

\begin{cor}
\label{teorema34}
For a positive integer $n$, suppose that $p^{a}$ divides~$n$ where $p$ is a prime and $a>0$. If $q$ is a prime such that $n/2<q<n$ and $n-q < p^{a}$, then $n$ satisfies Condition~$1$ with $p$ and~$q$.
\end{cor}

\begin{proof}
Pick $d=1$ in Theorem~\ref{MIAU}.
\end{proof}

Note that, under the assumptions of Corollary~\ref{teorema34}, the equality $n-q=p^a$ cannot hold, since $p$ divides $n$ and $p\ne q$ because $q$ does not divide~$n$. Hence there remains to study the case when $n-q>p^a$ and $q$ is the largest prime smaller than~$n$ while $p^a$ is the largest prime power dividing~$n$.
In other words, Condition~1 holds for $n$ whenever there is a prime between $n-p^a$ and~$n$.

The sequence of integers $n$ for which there is no prime between $n-p^a$ and $n$ can be found in The On-Line Encyclopedia of Integer Sequences (OEIS \cite{oeis}) with the reference A290203. Its first terms are the following:
\begin{equation}
\label{A290203}
126, 210, 330, 630, 1144, 1360, 2520, 2574, 2992, 3432, 3960, 4199 \dots
\end{equation}

\emph{Banderier's conjecture} \cite{banderier} claims that if $p_n\#$ denotes the $n$-th \emph{primorial}, that is, 
\[
p_n\#=p_1p_2\cdots p_n
\]
where $p_1,\dots,p_n$ are the first $n$ primes, and $q$ is the largest prime below $p_n\#$, then either $p_n\#-q=1$ or $p_n\#-q$ is a prime. 

\begin{proposition}
If Banderier's conjecture is true, then the sequence \eqref{A290203} contains all primorials $p_n\#$ such that $p_n\#-1$ is not a prime.
\end{proposition}

\begin{proof}
If $p_n\#-1$ is not a prime, then $p_n\#-q$ is a prime according to Banderier's conjecture. Since $p_n\#-q$ does not divide $p_n\#$, we infer that $p_n\#-q$ is bigger than~$p_n$, which is the largest prime power dividing $p_n\#$. 
\end{proof}

The first primorials $p_n\#$ such that $p_n\#-1$ is not a prime are 
\[
\text{$p_4\#=210$, \quad $p_7\#=510510$, \quad $p_8\#=9699690$, \quad $p_9\#=223092870$.}
\]
Inspecting this list could be a strategy to seek for a counterexample for Question~\ref{q}. 
The complementary list of primorials can be found in OEIS with reference A057704.

For any fixed value of~$d$, the number $\beta$ in Theorem~\ref{MIAU} is smallest when $q$ is as close as possible to~$n/d$. For this reason, we focus our attention on the largest prime $q_d$ below $n/d$ for various values of~$d$. This motivates the next definition.

\begin{definition1}
For positive integers $n$ and $1\le d<n/2$, let $q_d$ be the largest prime smaller than~$n/d$ and let $\beta_d=n-dq_d$. For each integer $c$ with $0\le c<(d+1)/2$, we call $[cq_{d},\, cq_{d} + \beta_d]$ a \textit{dangerous interval}.
\end{definition1}

By Theorem~\ref{MIAU}, if we attempt to prove that Condition~1 holds with $p$ and~$q_d$ assuming that $q_d>n/(d+1)$ ---that is, assuming that the dangerous intervals are disjoint--- we only need to care about values of $k$ that lie in a dangerous interval and are multiples of the largest power of $p$ dividing~$n$. 

In the case $d=1$, the only dangerous interval below $n/2$ is $[0,n-q_1]$.
When $d=2$, we have that $[0,n-2q_2]$ and $[q_2,n-q_2]$ are dangerous intervals. Since $n-q_2>n/2$, the second interval may be replaced by $[q_2,n/2]$ to carry our study further, as we do in the next section.

\begin{example}
The largest prime below $n=p_7\#=510510$ is $q_1=510481$ and the largest prime dividing $n$ is $p=17$. Here $n-q_1=29$ and therefore $\binom{n}{k}$ is divisible by $17$ or $510481$ for all $k$ except for $k=17$. 

On the other hand, the largest prime below $n/2=255255$ is $q_2=255253$. Thus $\beta_2=n-2q_2=4$ and therefore $[0,4]$ and $[255253,\,255257]$ are dangerous intervals. The second interval contains a multiple of $17$, namely $n/2$. However, since
\begin{align*}
    510510 & = 6\cdot 17^4 + 1\cdot 17^3 + 15\cdot 17^2 + 8\cdot 17\\
    255255 & = 3\cdot 17^4 + 0\cdot 17^3 + 16\cdot 17^2 + 4\cdot 17,
\end{align*}
we infer from Lucas' Theorem that $\binom{510510}{255255}$ is divisible by~$17$. Consequently, $\binom{n}{k}$ is divisible by $17$ or $255253$ for all~$k$. 
\end{example}

\section{Using the nearest prime below $n/2$}
\label{secq2}

Nagura showed in \cite{nagura} that, if $m \geq 25$, then there is a prime between $m$ and $(1 + 1/5)m$. Therefore, there is a prime $q$ such that $5n/6 < q < n$ when $n \geq 30$.
This implies that, if $n \geq 30$ and the largest prime-power divisor $p^a$ of $n$ satisfies $p^{a} \geq n/6$, then there is a prime $q$ between $n-p^a$ and~$n$ and hence Condition~1 holds for $n$ with $p$ and~$q$.

The following result is sharper.

\begin{prop}
If $n \geq 2010882$ and the largest prime-power divisor $p^a$ of $n$ satisfies $p^a\geq n/16598$, then $n$ satisfies Condition~$1$ with $p$ and the nearest prime $q$ below~$n$. 
\end{prop}

\begin{proof}
Schoenfeld proved in \cite{lowell} that for $m \geq 2010760$ there is a prime between $m$ and $(1 + 1/16597)m$. 
Hence, if $n \geq 2010882$ and the largest prime-power divisor $p^a$ of $n$ satisfies $p^{a} \geq n/16598$ then there is a prime between $n-p^a$ and~$n$, and therefore Condition~1 holds for $n$  by Corollary~\ref{teorema34}.
\end{proof}        

The following are consequences of Nagura's and Schoenfeld's bounds.

\begin{lemma}
\label{lemaxx}
Let $q_d$ be the largest prime below $n/d$ for positive integers $n$ and~$d$.
\begin{itemize}
\item[(a)]
If $n \geq 120$ and $d<5$, then $n/(d+1)<q_d$.
\item[(b)]
If $n \geq 3.34\cdot 10^{10}$ and $d < 16597$, then $n/(d+1)<q_d$.
\end{itemize}
\end{lemma}

\begin{proof}
By Nagura's bound \cite{nagura}, if $n/d \geq 30$, then $5n/6d < q_{d} < n/d$. Therefore, $n-dq_{d}<n/6$. If $d<5$, then $6d<5(d+1)$ and hence
\[
n < \dfrac{5n(d+1)}{6d}<q_{d}(d+1),
\] 
as claimed.
The proof of part~(b) is analogous using Schoenfeld's bound \cite{lowell}.
\end{proof}

In order to apply Theorem~\ref{MIAU} with $d=2$ for a given~$n$, we need that there is a prime $q$ such that $n/3<q<n/2$. If $q_2$ denotes the nearest prime below $n/2$, then the inequality $n/3<q_2$ holds if $n\geq 120$ by Lemma~\ref{lemaxx}. Since by \eqref{A290203} we have that $n-q_1<p^a$ if $n<126$, we may assume that $n/3<q_2$ without any loss of generality.

Note that the inequality $n/3<q$ is equivalent to $n-2q<q$, so the intervals $[0,n-2q]$ and $[q,n-q]$ are disjoint.

\begin{thm}
For an odd positive integer $n$ and a prime power $p^a$ dividing~$n$, suppose that there is a prime $q$ with $n/3<q<n/2$ and $n-2q < p^{a}$. Then $n$ satisfies Condition $1$ with $p$ and~$q$.
\end{thm}

\begin{proof}
By Theorem~\ref{MIAU}, in order to infer that $\binom{n}{k}$ is divisible by $p$ or~$q$,
the only cases that we need to discuss are those values of $k$ that are multiples of $p^a$ with $k\in[0,n-2q]$ or $k\in[q,n-q]$. By assumption, there are no multiples of $p^a$ in $[0,n-2q]$. Since $n-q>n/2$, we may focus on the interval $[q,n/2]$.
Since $n$ is odd, $n/2$ is not an integer; hence we are only left to prove that there is no multiple $k$ of $p^{a}$ with $q\leq k< n/2$. We will prove this by contradiction.

Thus suppose that $q \leq \lambda p^{a} < n/2$ for some integer~$\lambda$. The assumption that $n-2q<p^a$ implies that $n-p^a<2q$ and hence
\[
n/2-p^a/2<q\le\lambda p^a.
\]
Consequently,
$\lambda p^a < n/2 < (\lambda+1/2)p^a$.
If we now write $n=mp^a$, we obtain that $2\lambda<m<2\lambda+1$, which is impossible for an integer~$m$.
\end{proof}

The rest of this section is devoted to the case when $n$ is even. 

\begin{lemma}
\label{lema52}
Suppose that $n$ is even and there is a prime $q$ with $q<n/2$ and $n-2q < p^{a}$, where $p^a$ is the largest power of $p$ dividing~$n$.
If there is a multiple $k$ of $p^{a}$ in the interval $[q,n/2]$, then $p$ is odd and~$k=n/2$.
\end{lemma}

\begin{proof}
Suppose first that $p$ is odd. Then the integer $n/2$ is a multiple of $p^{a}$, so we may write $n/2=\lambda p^{a}$ for some integer~$\lambda$. If there is another multiple of $p^{a}$ in the interval $[q,n/2]$, then $q \leq (\lambda - 1) p^{a} < n/2$, and this implies that 
\[
n/2 - p^{a}=\lambda p^a-p^a= (\lambda-1)p^a \geq q. 
\]
Hence $n - 2q \geq 2p^{a}$, which is incompatible with our assumption that $n - 2q < p^{a}$.

In the case $p=2$ (so that $2^a$ is the largest power of $2$ dividing~$n$), 
we have that $n/2$ is divisible by $2^{a-1}$, and we may write $n/2=\lambda 2^{a-1}$ with $\lambda$ odd. If there is a multiple of $2^{a}$ in the interval $[q,n/2)$, then $q \leq \mu 2^{a} < n/2$, so $\mu<\lambda/2$ and $\mu\leq (\lambda-1)/2$ because $\lambda$ is odd. Therefore
\[
n/2 - 2^{a-1} =(\lambda-1)2^{a-1}\geq \mu 2^a\geq q. 
\]
Hence, as above, $n - 2q \geq 2^{a}$, which contradicts that $n - 2q < 2^{a}$.
\end{proof}

\begin{thm}
\label{central}
For an even positive integer~$n$, suppose that there is a prime $q$ with $n/3<q<n/2$ and $n-2q < p^{a}$, where $p^a$ is the largest power of $p$ dividing~$n$.
\begin{itemize}
\item[(a)] If $p=2$, then $n$ satisfies Condition~$1$ with $2$ and~$q$.

\vskip 0.1cm

\item[(b)] If $p\ne 2$, then $n$ satisfies Condition~$1$ with $p$ and~$q$ if and only if $\binom{n}{n/2}$ is divisible by~$p$.
\end{itemize}
\end{thm}

\begin{proof}
By Theorem~\ref{MIAU} and Lemma~\ref{lema52}, the only case left is $k=n/2$ for $p$ odd. Consequently, if $\binom{n}{n/2}$ is divisible by~$p$, then $n$ satisfies Condition~1 with $p$ and~$q$. Moreover, $\binom{n}{n/2}$ is not divisible by~$q$, since the base~$q$ expansions of $n$ and $n/2$ are, respectively, $2\cdot q + (n-2q)$ and $1\cdot q+(n/2-q)$. Hence the assumption that $\binom{n}{n/2}$ be divisible by~$p$ is necessary.
\end{proof}

Our last remarks in this section correspond to the case when $n$ is even, and they are only relevant if $p\ne 2$, by Theorem~\ref{central}. Sufficient conditions are given to infer that a prime $p$ divides $\binom{n}{n/2}$. The greatest integer less than or equal to a real number $x$ is denoted by $\lfloor x\rfloor$, and we write $v_{p}(n)=a$ if $p^a$ is the maximum power of $p$ such that $p^{a}$ divides~$n$.

Recall from \cite{legendre} that
\begin{equation}
\label{vpnfact}
v_{p}(n!) = \sum\limits_{k=1}^\infty \left\lfloor\dfrac{n}{p^{k}}\right\rfloor
= \dfrac{n-s_{p}(n)}{p-1},
\end{equation}
where $s_{p}(n)$ denotes the sum of all the digits in the base $p$ expansion of~$n$.

\begin{prop}
\label{prop44}
Suppose that $n$ is even.
A prime $p$ divides $\binom{n}{n/2}$ if and only if at least one of the numbers $\lfloor n/p^{r} \rfloor$ with $r \geq 1$ is odd.
\end{prop}

\begin{proof}
By comparing $v_{p}(n!)$ and $v_{p}((n/2)!)$ we see that, for each~$r$, 
\[
\left\lfloor\dfrac{n}{p^{r}}\right\rfloor = 2\left\lfloor\dfrac{n/2}{p^{r}}\right\rfloor
\] 
if $\lfloor n/p^{r } \rfloor$ is even. If $\lfloor n/p^{r} \rfloor$ is even for all~$r$, we conclude that $v_{p}(n!) = 2v_{p}((n/2)!)$, and hence $p$ does not divide $\binom{n}{n/2}$. However, if  $\lfloor n/p^{r} \rfloor$ is odd, then 
\[
\left\lfloor\dfrac{n}{p^{r}}\right\rfloor = 2\left\lfloor\dfrac{n/2}{p^{r}}\right\rfloor + 1
\] 
and consequently $v_{p}(n!)$ is greater than $2v_{p}((n/2)!)$.
\end{proof}

\begin{cor}
If $n$ is even and $(n-s_{p}(n))/(p-1)$ is odd, then $p$ divides $\binom{n}{n/2}$.
\end{cor}

\begin{proof}
This follows from Proposition~\ref{prop44} and Legendre's formula \eqref{vpnfact}.
\end{proof}

\begin{cor}
Suppose that $n$ is even.
\begin{itemize}
\item[(a)]
If any of the digits in the base $p$ expansion of $n/2$ is larger than $\left\lfloor p/2 \right\rfloor$, then $p$ divides $\binom{n}{n/2}$.
\item[(b)]
If one of the digits in the base $p$ expansion of $n$ is odd, then $p$ divides $\binom{n}{n/2}$.
\end{itemize}
\end{cor}

\begin{proof}
If a digit of $n/2$ in base $p$ is larger than $\left\lfloor p/2 \right\rfloor$, then when we add $n/2$ to itself in base~$p$ to obtain $n$ there is at least one carry. Similarly, if $n$ has an odd digit in base~$p$, then there is a carry when adding $n/2$ and $n/2$ in base~$p$. Hence, by Kummer's Theorem \cite{kummer} with $k=n/2$, if there is at least one carry when adding $n/2$ to itself in base~$p$, then $p$ divides $\binom{n}{n/2}$.
\end{proof}

\begin{cor}
Let $n$ be an even positive integer. Suppose that there is a prime~$q$ such that $n/3<q<n/2$ and $n-2q < p^{a}$, where $p^a$ denotes the largest power of $p$ dividing~$n$. If $p^{\left\lfloor\log n/\log p\right\rfloor} > n/2$,
then $p$ divides $\binom{n}{n/2}$ and therefore $n$ satisfies Condition~$1$ with $p$ and~$q$.
\end{cor}

\begin{proof}
The largest value of $r$ such that $p^{r} < n < p^{r + 1}$ is 
$\left\lfloor\log n/\log p\right\rfloor$.
Therefore, in Proposition~\ref{prop44}, the exponent $r$ is bounded by 
$\left\lfloor\log n/\log p\right\rfloor$.
Also note that $r \geq a$, where $a$ is the largest exponent of $p$ such that $p^a$ divides~$n$. If 
$p^{\left\lfloor\log n/\log p\right\rfloor}>n/2$,
then $\lfloor n/p^{r} \rfloor = 1$. Because this is odd, $p$ divides $\binom{n}{n/2}$ by Proposition~\ref{prop44}. 
\end{proof}

In those cases when the inequalities $n-q_1<p^a$ and $n-2q_2<p^a$ both fail for the largest prime power $p^a$ dividing~$n$, a possible strategy is to analyze the inequality $n-dq_{d}<p^a$ for bigger values of~$d$, where $q_{d}$ is the largest prime below~$n/d$.

Up to $1{,}000{,}000$ there are $88$ integers that do not satisfy $n-2q_{2}<p^a$, where $p^a$ is the largest prime power dividing~$n$. The On-Line Encyclopedia of Integer Sequences has published these numbers \cite{oeis2} with the reference A290290. Among these, there are $25$ that do not satisfy the inequality $n-3q_{3}<p^a$; there are $7$ that do not satisfy the inequality $n-4q_{4}<p^a$ either; there are $5$ for which the inequality $n-5q_{5}<p^a$ also fails, and there is only one integer for which the inequality $n-6q_{6}<p^a$ still fails (namely, $n=875160$). However, the value of $n-dq_d$ need not decrease as $d$ grows, and the number of dangerous intervals that one needs to inspect when $n-dq_d<p^a$ increases linearly with~$d$. 
Therefore this strategy is not conclusive, although it often works in practice.

\begin{example}
The largest prime power dividing $n=p_{14}\#=13082761331670030$ is $p=43$. In this case, $n-q_1=89$ and $n-2q_2=268$. Thus, Condition~1 fails for $p$ and $q_1$ and it also fails for $p$ and~$q_2$. Nevertheless, $n-3q_3=27$ works, as the dangerous interval $[q_3,n-2q_3]$ contains one multiple of~$43$, namely $n/3$, 
and $\binom{n}{n/3}$ is divisible by~$43$. Therefore Condition~1 holds for $p=43$ and~$q_3=4360920443890001$.
\end{example}

\begin{example}
For $n=210$, the inequality $n-q_1<7$ fails while $n-2q_2<7$ is true. However, 
$\binom{210}{105}$ is not divisible by~$7$. Hence we look for greater values of $d$ and find that $n-5q_5<7$ with $q_5=41$. Now $42\in[41,46]$ and $84\in[82,87]$, yet $\binom{210}{42}$ and $\binom{210}{84}$ are both divisible by~$7$. Hence Condition~1 is satisfied with $p=7$ and $q_5=41$.
\end{example}

\begin{example}
For $n=875160$, the inequality $n-dq_d<17$ is satisfied with $d=11$ but not with any smaller value of~$d$. There are $6$ dangerous intervals of length $n-11q_{11}=11$. Each of these intervals (except the first) contains one multiple of~$17$, and in each case the corresponding binomial coefficient $\binom{n}{k}$ happens to be divisible by~$17$. Therefore Condition~1 is satisfied with $p=17$ and $q_{11}=79559$.
\end{example}


\section{On the $N$-variation of Condition 1}
\label{sec4}
Recall from Definition~\ref{def1} that $n$ satisfies the \emph{$N$-variation} of Condition~1 if there are $N$ primes $p_1,\dots,p_N$ such that if $1 \leq k \leq n-1$ then $\binom{n}{k}$ is divisible by at least one of $p_1,\dots,p_N$.

\begin{thm}
If an even positive integer $n$ satisfies $n-2q<p^{a}$ for a prime $q$ with $n/3<q<n/2$, where $p^a$ is the largest power of $p$ dividing $n$ and $p \neq 2$, then $n$ satisfies the $3$-variation of Condition~$1$ with $p$, $q$ and any prime that divides $\binom{n}{n/2}$.
\end{thm}

\begin{proof}
According to part (b) of Theorem~\ref{central}, the only binomial coefficient $\binom{n}{k}$ with $1\leq k\leq n-1$ that might fail to be divisible by $p$ or $q$ is $\binom{n}{n/2}$. Hence it suffices to add an extra prime with this purpose.
\end{proof}

\begin{prop}
\label{p1p2}
For a positive integer $n$, let $q_1$ be the largest prime smaller than~$n$, let 
$p_{1}^{a_1}$ be the largest prime-power divisor of $n$ and let $p_{2}^{a_{2}}$ be the second largest prime-power divisor of~$n$. If $p_{1}^{a_{1}}p_{2}^{a_{2}} > n-q_1$, then $n$ satisfies the $3$-variation of Condition~$1$ with $p_{1}$, $p_{2}$ and~$q_1$.
\end{prop}

\begin{proof}
By Lucas' Theorem, for any $k$ such that $1 \leq k < p_{1}^{a_{1}}$, the binomial coefficient $\binom{n}{k}$ is divisible by~$p_{1}$, and for any $k$ such that $n-q_1 < k \leq n/2$ the binomial coefficient $\binom{n}{k}$ is divisible by~$q_1$. Thus we need to add a prime that divides at least the binomial coefficients $\binom{n}{k}$ with $p_{1}^{a_{1}} \leq k \leq n-q_1$ in which $k$ is a multiple of~$p_{1}^{a_{1}}$. 
For this, we pick~$p_2$ and therefore we only need to consider those values of $k$ that are, in addition, multiples of~$p_{2}^{a_{2}}$. 
The least $k$ that is a multiple of both prime powers is $p_{1}^{a_{1}}p_{2}^{a_{2}}$. 
Therefore, if $p_{1}^{a_{1}}p_{2}^{a_{2}} > n-q_1$, then all values of $k$ lying in the interval $p_{1}^{a_{1}} \leq k \leq n-q_1$ are such that $\binom{n}{k}$ is divisible by $p_{1}$ or~$p_{2}$. 
\end{proof}

In the statement of Proposition~\ref{p1p2}, the condition that $p_{1}^{a_{1}}p_{2}^{a_{2}} > n-q_1$ holds by Nagura's bound \cite{nagura} if we impose instead that $p_{1}^{a_{1}}p_{2}^{a_{2}} > n/6$.

For each~$n$, we are interested in the minimum number $N$ of primes such that $n$ satisfies the $N$-variation of Condition~1. We next discuss upper bounds for~$N$.

\begin{prop}
\label{d2}
For positive integers $n$ and~$d$, suppose that there is a prime $q$ such that $n/(d+1)<q<n/d$ and a prime-power divisor $p^a$ of $n$ such that $n-dq<p^a$. 
Then $n$ satisfies the $N$-variation of Condition $1$ with $N = 2+ \left\lfloor d/2 \right\rfloor$.
\end{prop}

\begin{proof}
By Theorem~\ref{MIAU}, the binomial coefficients $\binom{n}{k}$ are divisible by $q$ except possibly if $k$ lies in a dangerous interval. In the dangerous intervals we only need to consider those integers that are multiples of~$p^{a}$, since otherwise $\binom{n}{k}$ is divisible by~$p$.
Since we are assuming that $n-dq<p^a$, we know that in each dangerous interval there is at most one multiple of~$p^{a}$. This means that the worst case is the one in which there is a multiple of $p^{a}$ in every dangerous interval $[cq,cq+\beta]$ with $1\leq c \leq \left\lfloor d/2 \right\rfloor$. Hence we pick one extra prime for each such interval.
\end{proof}

\begin{cor}
If $1<d<5$ and $p^a > q_d+\beta_d$ where $p^a$ divides~$n$ and $q_d$ is the largest prime below~$n/d$, and $\beta_d=n-dq_d$,
then $n$ satisfies Condition~$1$ with $p$ and~$q_{d}$.
\end{cor}

\begin{proof}
By Lemma~\ref{lemaxx}, we may assume that $n/(d+1)<q_d$.
If $1<d < 5$, then $\left\lfloor d/2 \right\rfloor$ equals $1$ or~$2$. If $\left\lfloor d/2 \right\rfloor=1$, then the assumption that $p^{a} > q_{d} + \beta_d$ implies that no multiple of $p^{a}$ falls into any dangerous interval until~$n/2$. If $\left\lfloor d/2 \right\rfloor=2$, then we need to check that $2p^{a} > 2q_{d} + \beta_d$ in order to ensure that $2p^{a}$ does not fall into the third dangerous interval. The minimum value of $p^{a}$ such that our assumption $p^{a} > q_{d} + \beta_d$ holds is $q_{d} + \beta_d + 1$. The next multiple of $q_{d} + \beta_d + 1$ is $2q_{d} + 2\beta_d + 2$, which is greater than $2q_{d} + \beta_d$ and therefore $2p^{a}$ does not fall into the third dangerous interval.
\end{proof}

In order to refine the conclusion of Proposition~\ref{d2}, we consider the Diophantine equation 
\begin{equation}
\label{Diophantine}
p^{a}x - q_{d}y = \delta,
\end{equation}
for $0 \leq \delta \leq \beta_d=n-dq_d$, where $p^a$ is a prime-power divisor of a given number $n$ and $q_d$ is the largest prime below $n/d$ with $d\ge 1$. We keep assuming, as above, that $q_d>n/(d+1)$. We will also assume that $p\ne q_d$, which guarantees that \eqref{Diophantine} has infinitely many solutions for each value of~$\delta$. Specifically, if $(x_1,y_1)$ is a particular solution for some value of~$\delta$, then the general solution for this $\delta$ is
\[
x = x_{1} + rq_{d}, \qquad
y = y_{1} + rp^{a},
\]
where $r$ is any integer. In the next theorem we denote by $N(\delta)$ the number of solutions $(x,y)$ of \eqref{Diophantine} with $x>0$ and $0\le y\le \lfloor d/2\rfloor$ for each value of $\delta$ with $0\le\delta\le\beta_d$. Thus $N(\delta)=0$ precisely when \eqref{Diophantine} has no solution $(x,y)$ subject to these conditions.

\begin{thm}
\label{N}
For positive integers $n$ and~$d$, suppose that the largest prime $q_d$ below $n/d$ satisfies $q_d>n/(d+1)$, and let $\beta_d=n-dq_d$.
Let $p^a$ be a prime power dividing~$n$ with $p\ne q_d$. 
Then $n$ satisfies the $N$-variation of Condition~$1$ with
\[
N = 2 + \sum\limits_{\delta = 0}^{\beta_d} N(\delta),
\]
where $N(\delta)$ is the number of solutions $(x,y)$ of $p^ax-q_dy=\delta$ with $x>0$ and $0\le y\le \lfloor d/2\rfloor$ for each value of $\delta$ with $0\le\delta\le\beta_d$.
\end{thm}

\begin{proof} 
The number $N(\delta)$ counts how many times a multiple of $p^a$ falls into a dangerous interval $[cq_d,cq_d+\beta_d]$ at a distance $\delta$ from the origin of that interval. Thus we pick an extra prime for each such case, and add two to the sum in order to account for $p$ and~$q_d$.
\end{proof}

\begin{example}
\label{darrer}
The largest prime-power divisor of $n=96135$ is $p=29$. For $d=4$ we find that $q_4=24029$ and $\beta_4=19$. Since $24029\equiv 17$ mod~$29$, the only solution $(x,y)$ of the Diophantine equation $29x-24029y=\delta$ with $x>0$ and $0\le y\le 2$ is $(829,1)$ for $\delta=12$. Thus, $N(12)=1$ and $N=3$ for $d=4$. In other words, the only occurrence of a multiple of $29$ in a dangerous interval for $d=4$ is $24041\in[24029,\,24048]$. This example shows that the bound $2+\lfloor d/2\rfloor$ given in Proposition~\ref{d2} can be lowered.
\end{example}

The number $N$ given by Theorem~\ref{N} is not a sharp bound. For those multiples $p^ax$ of $p^a$ falling into a dangerous interval $[cq_d,cq_d+\beta_d]$, it often happens that the corresponding binomial coefficient $\binom{n}{p^ax}$ is divisible by~$p$, as in Example~\ref{darrer} or in other examples given in the previous sections. It could also be divisible by $q_d$ if $d\ge q_d$. When $d<q_d$, we have that $n$ satisfies Condition~1 with $p$ and $q_d$ if and only if the binomial coefficient $\binom{n}{p^ax}$ is divisible by~$p$ for every solution $(x,y)$ of \eqref{Diophantine} with $x>0$ and $0\le y\le\lfloor d/2\rfloor$, since $n=dq_d+\beta_d$ and $p^ax=yq_d+\delta$ with $\delta\le\beta_d<q_d$ and $y\le\lfloor d/2\rfloor<d$, so $\binom{n}{p^ax}$ is not divisible by~$q_d$ by Lucas' Theorem. Note also, for practical purposes, that $\binom{n}{p^ax}\equiv\binom{n/p^a}{x}$ mod~$p$.

%

\section{Every number has multiples for which Condition~1 holds}
\label{sec5}

We next prove that every positive integer $n$ has infinitely many multiples for which Condition~1 holds. We are indebted to R.\,Woodroofe for simplifying and improving our earlier statement of this result, which was based on prime gap conjectures.

It follows from the Prime Number Theorem \cite{hl} that, given any real number $\varepsilon >0$, there is a prime between $m$ and $m(1+\varepsilon)$ for sufficiently large~$m$.
This fact can be used to prove the following:

\begin{thm}
For every positive integer $n$ and every prime~$p$,
the number $np^k$ satisfies Condition~$1$ with $p$ and another prime, for all sufficiently large values of~$k$.
\end{thm}

\begin{proof}
For any prime $p$ and any $k>0$, let $m=np^k-p^k=p^k(n-1)$. Then
\[
np^k=m+p^k=m\left(1+\frac{1}{n-1}\right).
\]
Therefore, by the Prime Number Theorem, there is a prime between $m$ and $np^k$ for all sufficiently large values of~$k$. Choose the largest prime $q$ with this property. Thus,
\[
np^k-p^k<q<np^k,
\]
so $np^k-q<p^k$, from which it follows, according to Corollary~\ref{teorema34}, that $np^k$ satisfies Condition~1 with $p$ and~$q$.
\end{proof}

\begin{thm}
For every positive integer $n$ there is a number $M$ such that if $p$ is any prime with $p>M$ then $np$ satisfies Condition~$1$ with $p$ and another prime.
\end{thm}

\begin{proof}
Given $n$, let $\varepsilon=1/(n-1)$. Choose $m_0$ such that there is a prime between $m$ and $m(1+\varepsilon)$ for all $m\geq m_0$, and let $M=\varepsilon m_0$. If $p$ is any prime such that $p>M$, then for $m=p(n-1)$ we have
\[
np=m+p=m\left(1+\frac{p}{m}\right)=m\left(1+\frac{1}{n-1}\right)=m(1+\varepsilon).
\]
Therefore, by our choice of $m_0$, there is a prime between $m$ and $np$. If $q$ is the largest prime with this property, then
$np-p<q<np$, and consequently $np$ satisfies Condition~1 with $p$ and~$q$.
\end{proof}

Prime gap conjectures provide information relevant to our problem. For example, if $p_i$ denotes the $i$-th prime, then Cram\'er's conjecture \cite{cramer} claims that
there exist constants $M$ and $N$ such that if ${p}_{i} \geq N$ then 
\[
{p}_{i+1} - {p}_{i} \leq M(\log {p}_{i})^{2}.
\]

\begin{prop}
\label{mfactors}
Let $m$ be the number of distinct prime factors of~$n$. If Cram\'er's conjecture is true and $n$ grows sufficiently large keeping $m$ fixed, then $n$ satisfies Condition~$1$.
\end{prop}

\begin{proof}
If $n$ has $m$ distinct prime factors, then $\sqrt[m]{n} \leq p^{a}$, where $p^{a}$ is the largest prime-power divisor of~$n$. 
Let $M$ and $N$ be the constants given by Cram\'er's conjecture. Pick $n_0$ such that if $n\geq n_0$ then $M(\log n)^{2} < \sqrt[m]{n}$. For every $n$ such that $n\geq n_0$ and $N\le p_i<n\leq p_{i+1}$ (where $p_i$ denotes the $i$-th prime), we have
\[
n-p_i\leq p_{i+1}-p_i\le M(\log p_i)^2<M(\log n)^2<\sqrt[m]{n}\leq p^a,
\]
from which it follows that $n$ satisfies Condition~1 with $p$ and $p_i$.
\end{proof}

We note that the argument used in the proof of Proposition~\ref{mfactors} yields an alternative proof of the fact that Condition~1 holds for a set of integers of asymptotic density~$1$ if Cram\'er's conjecture holds, a result first found in~\cite[\S\,5]{sw}:

\begin{thm}[\cite{sw}]
\label{density}
If Cram\'er's conjecture is true, then the set of numbers in the sequence \eqref{A290203} has asymptotic density zero.
\end{thm}

\begin{proof}
Suppose that Cram\'er's conjecture holds with constants $M$ and~$N$, and denote by $\omega(n)$ the number of distinct prime divisors of~$n$. Thus $n^{1/\omega(n)}\leq p^a$, where $p^a$ is the largest prime-power divisor of~$n$. According to \cite[\S\,3.2]{hr}, for every $\varepsilon>0$ the inequality
\begin{equation}
\label{hardy}
\omega(n)<(1+\varepsilon)\log\log n
\end{equation}
holds for all $n$ except those of a set of asymptotic density zero.
Since
\[
\lim_{n\to\infty}\frac{n^{1/\log\log n}}{(\log n)^k}=\infty
\]
for all~$k$, there is an $n_0$ such that $n^{1/\omega(n)}>M(\log n)^2$ if $n\ge n_0$. Now, if $n$ is bigger than $n_0$ and satisfies $N\le p_i<n\le p_{i+1}$, and moreover $n$ is not in the set of integers for which \eqref{hardy} fails, then
\[
n-p_i\leq p_{i+1}-p_i\le M(\log p_i)^2<M(\log n)^2<n^{1/w(n)}\leq p^a.
\]
Therefore, $n$ satisfies Condition~1 with $p$ and $p_i$.
\end{proof}

\section{Multinomials}
\label{sec6}

We also consider a generalization of Condition~1 to multinomials. We say that a positive integer $n$ satisfies \emph{Condition~$1$ for multinomials of order~$m$} if there are primes $p$ and $q$ such that the multinomial coefficient
\[
{n \choose k_{1},k_{2},\dots,k_{m}} = \frac{n!}{k_{1}!k_{2}! \cdots k_{m}!}
\]
is divisible by either $p$ or~$q$ whenever $k_{1}+ \dots + k_{m} =n$ with $1 \leq k_{i} \leq n-1$ for all~$i$. 

\begin{prop}
If $n$ satisfies Condition~$1$ with two primes $p$ and~$q$, then $n$ satisfies Condition~$1$ for multinomials of any order $m\le n$ with $p$ and~$q$.
\end{prop}

\begin{proof}
This follows from the equality
\[
{n \choose k_{1},k_{2},\dotsc,k_{m}} =
\binom{n}{k_{1}} \binom{n-k_{1}}{k_{2}} \binom{n-k_{1}-k_{2}}{k_{3}} \dotsm \binom{k_{m}}{k_{m}},
\]
and the fact that ${n \choose k_1}$ is divisible by $p$ or $q$ by assumption.
\end{proof}

Therefore, if Condition 1 is proven for binomial coefficients, then it automatically holds for multinomial coefficients.

\subsubsection*{Acknowledgements}
I am indebted to my mentor Oscar Mickelin for his guidance throughout this research and to Prof.\ Russ Woodroofe for correspondence and kind suggestions. This work was carried out during the Research Science Program at MIT in the summer of 2017 and was supported by the Center for Excellence in Education, the MIT Mathematics Department, and the Youth and Science Program of Fundaci\'o Catalunya La Pedrera (Barcelona).


\begin{thebibliography}{99}

\bibitem{banderier}
\url{https://lipn.univ-paris13.fr/~banderier/Computations/prime_factorial.html}
(last consulted on 11 August 2018).
\bibitem{bertrand} J. Bertrand, M\'emoire sur le nombre de valeurs que peut prendre une fonction quand on y permute les lettres qu'elle renferme, \textit{Journal de l'\'Ecole Royale Polytechnique}, 30:123--140, 1845.
\bibitem{oeis} S. Casacuberta, Sequence A290203 in The On-Line Encyclopedia of Integer Sequences, published electronically at \url{http://oeis.org/A290203}, 24 July 2017.
\bibitem{oeis2} S. Casacuberta, Sequence A290290 in the On-Line Encyclopedia of Integer Sequences, published electronically at \url{http://oeis.org/A290290}, 26 July 2017.
\bibitem{glucas} K. S. Davis and W. A. Webb, Lucas' theorem for prime powers, \textit{European Journal of Combinatorics}, 11:229--233, 1990.
\bibitem{cramer} A. Granville, Harald Cram\'er and the distribution of prime numbers, \textit{Scandinavian Actuarial Journal}, 1995:337--360, 1995.
\bibitem{hl} G. H. Hardy and J. E. Littlewood, Contributions to the theory of the Riemann zeta-function and the theory of the distribution of primes, \textit{Acta Mathematica}, 41:119--196, 1916.
\bibitem{hr}
G. H. Hardy and S. Ramanujan, The normal number of prime factors of a number $n$, \textit{Quart. J. Math.}, 48:76--92, 1917.
\bibitem{distribution} A. E. Ingham, \textit{The Distribution of Prime Numbers}, Cambridge University Press, 1932.
\bibitem{kummer} E. Kummer, {\"U}ber die Erg{\"a}nzungss{\"a}tze zu den allgemeinen Reziprozit{\"a}ts\-geset\-zen, \textit{Journal f{\"u}r die reine und angewandte Mathematik}, 44:93--146, 1852.
\bibitem{legendre} A. M. Legendre, \textit{Th\'eorie des Nombres}, Paris: Firmin Didot Fr\`eres, 1932.
\bibitem{lucas} \'E. Lucas, Th\'eorie des fonctions num\'eriques simplement p\'eriodiques, \textit{American Journal of Mathematics}, 44:184--196, 1878.
\bibitem{nagura} J. Nagura, \textit{On the interval containing at least one prime number}, Proceedings of the Japan Academy, 1952.
\bibitem{lowell} L. Schoenfeld, \textit{Sharper bounds for Chebyshev functions $\psi(x)$ and $\theta(x)$}, ii, Mathematics of Computation, 1976. 
\bibitem{sw} J. Shareshian and R. Woodroofe, Divisibility of binomial coefficients and generation of alternating groups, \textit{Pacific Journal of Mathematics}, 292:223--238, 2018.

\end{thebibliography}
\end{document}